# The Fate of Russian Translations of Cantor

Galina Sinkevich, St. Petersburg State University of Architecture and Civil Engineering, 4, Vtoraja Krasnoarmejskaja ul., St. Petersburg 190005, Russia

galina.sinkevich@gmail.com

Home address: Galina Sinkevich, P.O. 197101 box 60, St. Petersburg, Russia

**Biographical sketch.** Galina Sinkevich is associate professor of the Department of Mathematics, St. Petersburg State University of Architecture and Civil Engineering, Russia. She has been working on the history of mathematics for more than 30 years. She has research works in the history of the set theory and theory of functions; they include 3 monographs and 50 articles about Cantor, Sierpiński, Lusin, Weierstrass, Heine, Hankel, Dedekind, Dini, Stiefel, Pascal, Méray. Thanks to her efforts, the house where Cantor was born was found and a memorial plaque was installed. She conducts a research seminar on the History of Mathematics at the St. Petersburg Department of V.A. Steklov Institute of Mathematics of the Russian Academy of Sciences. Seminar site http://www.mathnet.ru/php/conference.phtml?option_lang=rus&eventID=10&confid=504

## Key words: set theory, Cantor, Russian translations.

**Highlights**:
1. All Russian translations of Cantor's systematized.
2. The authors of translations characterized with their works.
3. The unique material about cancelled publication Cantor in Russian founded out, it is its first communication.

**Abstract**. This is the history of translating Cantor's works into Russian from 1892 to 1985 in Odessa, Moscow, Tomsk, Kazan, S.-Petersburg, Leningrad. Mathematicians and philosophers in Russia took the ideas of the theory of sets enthusiastically. Such renowned scholars and scientists as Timchenko, Shatunovsky, Vasiliev, Florensky, Mlodzeevsky, Nekrasov, Zhegalkin, Yushkevich Sr., Fet, Yushkevich Jr., Kolmogorov, and Medvedev took part in their popularisation. In 1970 Academician Pontryagin rated the theory of sets as useless for young mathematicians, and the translated works of Cantor were not published. This article first describes the tragic fate of this translation.

From 1872 to 1897 Cantor wrote his basic works devoted to the theory of sets. Russian mathematicians who visited universities of Berlin and Gottingen and read Crelle's Journal, Mathematische Annalen, Acta Mathematica got to know the ideas of the theory of sets. Cantor's ideas gradually permeated research activities and teaching, appeared in press in the form of expositions and translations. We are going to review the history of Cantor's heritage in Russia from 1892 to 1985.

## Odessa. 1892. I.Y. Timchenko

We found the first references (1892) to Cantor's works in Russia in works of Ivan Y. Timchenko (1863–1939) who graduated from Novorossiysk University in Odessa in 1885 and subsequently became a professor in Odessa. Timchenko studied astronomy, mathematics and history of mathematics, travelled abroad to work in libraries (in 1890, 1892, 1893, and 1896). Timchenko chose historical analysis of development of the theory of analytical functions as the subject of his MPhil. His work entitled "Basis of the theory of analytical functions" was published in three editions of "Proceedings of the Department of Mathematics of Novorossiysk Scientists" in 1892 and 1899, and presented in 1899 [Timchenko 1892,1899].

His in-depth research covers the period from the ancient world to late 19$^{th}$ century. He considered development of basic ideas underlying the theory of analytical functions in this work. The most important of these ideas is the concept of continuity and related concepts of neighbourhood and limit point. Timchenko pays tribute to Weierstrass in the development of the concept of neighbourhood and uniform convergence of series, and to Georg Cantor in the geometrical treatment of the concept of continuity in his works devoted to linear manifolds. Timchenko points out the relationship between Cantor's apprehension of continuity ("integrity" according to Timchenko) and Leibniz' principle of continuity [Timchenko 1892, p. 12]. Notably, Timchenko addressed the key works of Cantor. The first one was his work of 1872, "Ueber die Ausdehnung eines Satzes aus der Theorie der trigonometrischen Reihen", where a new concept of a number and a concept of an limiting point were introduced. Mathematicians like G. Schwartz and W. Dini

[Dini, 1878] who used to give the course of analysis happily picked up this idea. The second and most famous one is the Fifth Memoir (Ueber unendliche lineare Punktmannigfaltigkeiten) of the cycle which comprises 6 parts published in 1879-1884. This work is entitled "Grundlagen einer allgemeinen Mannigfaltigkeitslehre". It contains all basic definitions and theorems, including the definition of the empty set, perfect set, concept of the real number and its benchmarking compared to similar concepts of Weierstrass and Dedekind; it introduces a power scale and sets the continuum hypothesis.

### Odessa. 1896. S.O. Shatunovsky

Samuil Osipovich Shatunovsky (1859–1929) was born to a poor family of a Jewish artisan. He did not obtain a systematic education, having attended various learning institutions in Russia and Switzerland. For a long time Shatunovsky earned his living by giving lessons in small provincial towns. He used to write works in math devoted to issues in geometry and algebra, axiomatic determination of a value; published them in Russian and foreign magazines, and had an eye on European mathematical life. Shatunovsky appeared in Odessa between 1891 and1893. Owing to Odessa professors, Shatunovsky was as an exception allowed taking examinations for Master's degree which opened him the door to teaching [Chebotarev, 1940]. Shatunovsky was only granted the position of privat-docent when he was almost 47, and after 1917 he became a professor of Novorossiysk University. Since 1905, Shatunovsky was reading analysis using definitions and methods of the theory of sets. The first concepts of mathematical analysis are set forth from the perspective of the theory of sets. However, Russian terms differed from the contemporary ones. Thus, a set, for example, is called a 'complex'. He constructed a real numbers system, introduced the concept of a convergent complex (dense set). He first introduced the concept of a removable discontinuity of a function. His terms were unique; however, his presentation was rigorous. The lectures were lithographed in 1906–1907 and republished in 1923 [Shatunovsky,

1923]. Among his students were G.M. Fichtenholz, D.A. Kryzhanovsky, and I.V. Arnold. One can unquestioningly see the effect of this course in "Fundamentals of Analysis" by G.M. Fihtenholz (1888–1959) who attended Shatunovsky's lectures and graduated from Novorossiysky University in Odessa in 1911.

Shatunovsky has demonstrated an amazing scientific undersense choosing works to be translated. He was the first to translate Dedekind's "Stetigkeit und irrationale Zahlen" (written by Dedekind in 1872 and translated by Shatunovsky in 1894) [Dedekind, 1894] and Cantor's "Ueber eine Eigenschaft des Inbegriffes aller reellen algebraischen Zahlen" (written by Cantor in 1874 and translated by Shatunovsky in 1896) [Shatunovsky, 1896]. It was in these two works that a new concept of a number was created to form basis for the $20^{th}$ century mathematics.

From 1886 to 1917, there was a journal in Odessa, "Bulletin of Experimental Physics and Elementary Mathematics"(Vestnik opytnoy fiziki i elementarnoy matematiki). In 1896, they published an article of O.S. Shatunovsky in issue 233 entitled "Proof of Existence of Transcendental Numbers (as per Cantor)" [Shatunovsky, 1896]. Shatunovsky had laid down the proof of theorems from Cantor's work of 1874, "Ueber eine Eigenschaft des Inbegriffes aller reellen algebraischen Zahlen", having added a description of his more recent achievements, in particular, the concept of power of set which appeared only in 1878 in Cantor's work "Ein Beitrag zur Mannigfaltigkeitslehre".

## Moscow. 1900. B.K. Mlodzeevsky

Owing to the acquisition of literature and research travel, Moscow mathematicians were aware of West-European scientific achievements. To obtain the Master's Degree, students used to attend lectures at German and French research centres for at least one term. Lectures given in Moscow University included information on scientific achievements. The Theory of Real Variable Functions was read by

Boleslav Kornelievich Mlodzeevsky (1858–1923). Owing to the theory of sets, the course of Mathematics and Theory of Functions in the first place was rearranged on other basis. Mlodzeevsky used the course of Ulisse Dini as a base. Ulisse Dini used Cantor's results in his course as early as in 1870s [Dini, 1878]. Mlodzeevsky gave his course in the autumn term of 1900 and thereafter gave it a couple of times until 1908 [Medvedev, 1986].

Notes of Mlodzeevsky's lectures delivered in 1902 were found in archives of P.A. Florensky, then third-year student of the Department of Mathematics. The course included 29 lectures (three lectures a week). Based on the lecture notes, F.A. Medvedev supposed that "Mlodzeevsky did not seem to be directly familiar with G. Cantor's works at that time. The name and numerous set-theoretic and function-theoretic results of the latter are repeatedly mentioned in the lectures. However, judging from the nature of these mentioning (lack of direct references to Cantor's works or clarifications to the effect that certain considerations were set forth based on one of the above-listed works, etc.), it is reasonable to suppose that by 1902, B.K. Mlodzeevsky had learnt of Cantor's work second-hand, mainly from works of P. Tanneri, G. Tanneri, and A. Schoenflies" [8, c. 134]. The theory of sets is used in Mlodzeevsky's lectures to present the theory of function argument. He considered point sets ("clusters of points") and functions on the sets; he introduced the concept of a limiting point and cluster set; classified sets into first and second kind; formulated a theorem stating that a measure of a set of the first kind equals zero; the upper and lower limit, the concept of power of sets, countability    of rational and polynomial numbers; equipotency of various dimentions continua, denumerability of a countable sum of countable sets; uncountability of continuum with reference to the continuum hypotheses; perfect sets; ordinal type ("specie"); well-order (a "well-organized group"); transfinite numbers and alephs [Medvedev, 1986, p. 138–139].

**Moscow. 1904. P.A. Florensky**

Pavel Alexandrovich Florensky (1882–1937), a prominent philosopher theologian and priest later shot dead, was from 1900 to 1904 a student of the Department of Mathematics of Moscow University. In the autumn term of academic year1902/03, he attended a course of Mlodzeevsky's lectures where he learnt of Cantor's theory of sets. Since 1903, Florensky was working on his thesis entitled "The Idea of Discontinuity as an Element of Outlook" a preamble to which was published in 1986 in Historical Studies in Math (Istoriko-matematicheskie issledovania) [Florensky, 1903]. Florensky writes about Cantor's representation of continuity.

Florensky turned to Cantor's theory for a second time in 1904 in his work entitled "Symbols of Eternity (A Sketch of Cantor's Ideas)" [Florensky, 1904]. Florensky made it his crusade to paraphrase the meaning of Cantor's works. He described the development of definitions of potential and actual continuity in history of philosophy and continued with description of Cantor's theory of cardinals. However, he basically addressed Cantor's most recent works he wrote to provide a philosophical underpinning of his understanding of the continuity and theory of kinds: "Ueber die verschiedenen Standpuncte in bezug auf das Actuelle Unendliche" (1886) and "Mitteilungen zur Lehre vom Transfiniten" (1888). Cantor's theory, although in a concise form, was set forth appropriately. Basically, it is a paraphrase of the two articles named above. Florensky focused on the philosophical aspect of the theory more inclined to philosophy of religion. He gives Cantor credit for the introduction of actual infinity symbols. Further Florensky tries to understand Cantor's scientific motivation the background whereof he searches for in Cantor's biography, although he himself admits that "Cantor's biographic information has never been published and therefore, facts are extremely scarce. So one has to interpolate intuitively. However, having created a vision of Cantor's personality *of one's own*, it becomes extremely difficult to prove the rightfulness of one's own vision." [Florensky, 1904, p. 120, emphasis in original]

Florensky attributes the determination and purposeful nature of Cantor's scientific track to Jewish religiosity enhanced to self-sacrifice. We may make an allowance for the young age of Florensky who was but twenty-two when he started interpolating or rather imaging his views (and those of Vl. Soloviev) on the inner world of a scientist he did not know confusing the ideas of ethnicity and religious affiliation. Now we already know that Cantor was Lutheran who was born to a family of Lutheran father and Catholic mother; that it was merely his agnate grandfather who was Jewish, and in the next generation his father, brother and sister were Lutheran, and one of his father's other sisters was Orthodox Christian. His male line goes up to Portuguese Jews who settled in Copenhagen; his female line goes up to Austrian Czechs and the Hungarians who were Catholic [Sinkevich, 2012, 2014]. With parents belonging to different confessions Georg Cantor was not very religious and later, he consulted only Catholic theologians in search for a theological substantiation for concepts of his theory, although appealed on a point of all philosophic literature devoted to issues of eternity and continuum.

Florensky believed that one can comprehend the reality in symbols [Florensky, 1904, p. 126] and therefore absolutized Cantor's aspiration to create transfinite symbols. "Whereas Cantor as an individual is a most real-life image of a Jew, his philosophy is pretty much the same." [Florensky, 1904, p. 127]. Confusing ethnic and religious characteristics again, Florensky concludes that this was the reason why Cantor considered the actual infinity: "The idea of a complete infinity for both the absolute individual, the God, and for a human is the domain of Jewry, and this idea seems to be the most material grounds for Cantor. Meanwhile, others, the Aryans, acknowledge only potential infinity, "bad", indefinite and infinite, the very thought of nonexistence of the actual infinity seems unmerciful to his sole." [Florensky, 1904, p. 127].

Having graduated from the University, Florensky entered the Moscow Ecclesiastical Academy to become a priest.

By the first decade of the 20<sup>th</sup> century, Cantor's theory had spread in mathematic community of Europe and Russia. On the basis thereof, the theory of measure originated in works of Borel, Lebesgue, and Baire. In 1911, the school of the theory of functions and thereafter, the school of descriptive set theory started to form in Moscow. D.F. Egorov and N.N. Luzin were among their originators.

**Kazan. 1904–1908. A.V. Vasiliev.**

Since 1874, having graduated from Petersburg University, Alexandre Vasilievich Vasiliev (1853–1929) worked at Kazan University first as a privat-docent and from 1887, as a professor. His broad education, proficiency in languages, and numerous contacts with foreign scientists enabled him to make a good organizer and enlightener. He engaged in both research and socio-political activity, and advocated Lobachevsky's ideas having prepared his collected works edition for publication. Georg Cantor's uncle, Dmitry Ivanovich Meyer (1819–1856), famous lawyer and creator of Russian civil law [Sinkevich, 2012], worked in Kazan until 1855. There were two portraits in Vasiliev's study: a portrait of Lobachevsky and that of Meyer. Vasiliev knew Cantor from letters they exchanged and advocated his ideas.

From 1904 to 1908 A.V. Vasiliev's "Introduction into Analysis" was published in the publishing office of Kazan University to set forth principles of the theory of sets. According to S.S. Demidov, "Little by little, courses of analysis started to shape into the present-day courses of the kind. Mathematicians from Odessa (S.O. Shatunovskiy), Kiev (B.Y. Bukreev), and Kazan (A.V. Vasiliev) were the first to do so" [Demidov, 2000, p.77].

Cantor's theory was about 30 years of age by 1907. Over this time, works of his successors and criticism of opponents enriched the theory and thus finally shaped it. However, 10 basic articles of Cantor made the entire theory. The first summary monograph of Schoenflis appeared in 1900. However, this one was not complete either. The theory had to be presented in its entirety.

Cantor's theory of sets includes two parts: the theory of linear point sets and theory of transfinite numbers. Mlodzeevsky set the task of comprehensive presentation of the theory to two of his candidates: V.L. Nekrasov and I.I. Zhegalkin. Nekrasov was supposed to present the theory of point sets in detail and Zhegalkin, the theory of transfinite numbers. Each had to add results of his own to the presentation. Both candidates met the challenges.

Zhegalkin defended his thesis on 12 March, and Nekrasov on 4 October. Their Master's dissertations were published a year earlier to become the Russia's first monographs in the theory of sets.

## Moscow-Tomsk. 1907. V.L. Nekrasov

Vladimir Leonidovich Nekrasov (1864–1922) graduated from Kazan University where he stayed to work as a teacher. However, in 1900 he was transferred to the newly formed Tomsk Institute of Technology, Department of Abstract Mathematics. In order to prepare the Master's dissertation, in 1902-1903 he stayed in Europe on an academic mission. His Master's dissertation entitled "The Geometry and Measure of Linear Point Domains[1]" was published in 1907 in «News-Bulletin of Tomsk Institute of Technology" (Izvestija Tomskogo Technologicheskogo Instituta) [Nekrasov, 1907].

Chapter 1 contains a detailed historical sketch of basic results of the theory of sets and theory of measure, and an exhaustive bibliographical review. Before "Einleitung in die Mengenlehre" by A. Fraenkel appeared in 1919, Nekrasov's bibliography was the most complete. The list of references arranged chronologically from 1638 to 1907. In the third chapter entitled "The Most Recent Works", Nekrasov supplemented this list with references to works which had appeared by the time the manuscript of the third chapter went to the press. Nekrasov was trying to separate the theory of point sets from that of abstract sets: "As far as the size is concerned, already Cantor found that point ranges may be finite, countable, or have the size of continuity. Finding out relation of the latter size within the series of *alephs* is the task of the theory of transfinite numbers which is none of our interest here" [ibid., p. 98]. Starting to review history from discovery of infinitely small by Newton and Leibniz, Nekrasov writes that "*Bolzano* was father of contemporary theory of domains. However, it was *G. Cantor* who developed and made it strictly scientific." [ibid., p. 2, emphasis as in the original]. In the third chapter, Nekrasov supplements the number with precursors and Galilei with his example showing the correspondence of infinite sets of a natural number and their squares [ibid., p. 225]. Nekrasov lays emphasis on concepts of a limiting point and arbitrary set introduced by Cantor as fundamental. Further, he gives basic provisions of the theory of point sets and names the three basic characteristics of linear domains: size, structure, and measure.

---

[1] Nekrasov's "domain" shall be understood as a set.

The second chapter contains Nekrasov's own results related to the geometry of linear sets which correspond to the three types of deployment and combinations thereof for both closed and open sets. The structure of point of discontinuity of functions is Nekrasov's applied results. We can't but mention that Nekrasov was perhaps the first to note Ulisse Dini's priority in classification of the points of discontinuity [ibid, p. 102]. The third chapter was supplemented with new references and historical ordering of the development of the set theory ideas. The fourth chapter provides measure theory of A. Lebesgue and W. Young, although Nekrasov takes out the beginning of the measure theory from Riemann and Hankel. Nekrasov noted that Cantor's theory was recognized: "The right to exist and the role of the theory of domains within the general system of science has been established: this theory is reckoned and now its effect cannot be avoided quite in a number of branches of analysis. And this entire evolution happened within some 30 years, let alone so-to-say its prehistoric period." [ibid., p. 97].

Owing to the thorough historical analysis, elaborate presentation of Cantor's theory of point sets, and Nekrasov's own results, the monograph remains significant to the present day.

## Moscow. 1907. I.I. Zhegalkin

Having graduated from Moscow University, in 1906-1907 Ivan Ivanovich Zhegalkin (1869–1947) held a course in abstract set theory; in 1907 published his monograph entitled "Transfinite Numbers" [Zhegalkin, 1907], and in 1908 defended his Master's dissertation with the same title. Thereafter he headed research in mathematical logic where he obtained substantial results, having connected classic logic and residue arithmetic modulo 2. Residue ring modulo 2 is used to be called Zhegalkin algebra. In his subsequent works he proved monadic predicate calculus solvability.

Zhegalkin presents Cantor's algebra of transfinite numbers in his dissertation in his own way, deductively. There is no list of references, but for a couple of references to works of Cantor, Dedekind, Zermelo, and Bernstein. He mainly provides a transformed presentation of Cantor's last article of 1897 entitled "In support of the theory of transfinite sets". Zhegakin's presentation of the preamble is different. He thus hopes to avoid those contradictions in the theory which have come to light by the start of the $20^{th}$ century and were associated with the problem of well-ordering and Zermelo theorem. Zhegalkin adds more stringent arguments to Cantor's proof. He deserves credit for the statement regarding independence of the problem of choice from all other mathematic axioms made by him long before Serpinsky and Goedel.

In the first chapter, Zhegalkin tried to build a theory of cardinal and ordinal numbers before introducing the concept of finite and infinite. He introduces the concept of a finite set, ordering, and well-ordering; concept of a sum, product, and mapping of sets; and his own concept of a "not genuine" set. In the second chapter, he considers relation of equivalence, power (as a cardinal number), addition, multiplication operations, and raising to power, and this was the end of the theory of powers. The third and fourth chapters are devoted to the concept of an ordered set and concept of a type, and their respective properties. The fifth chapter considers a completely ordered set and Zermelo theorem ("every set can be thought as well-ordered set"). Zhegalkin proves the possibility to order the family of sets for the case of disjoint sets. (Zhegalkin calls them "detached").

The sixth chapter studies properties of ordinal numbers, i.e. types of completely ordered sets. Only after the theory of cardinal and ordinal numbers has been built, he considers finite sets and numbers in the seventh chapter; in the eighth chapter, he extracts countable sets therefrom as sets of all finite numbers. In the ninth chapter, he introduces congruence of powers; Chapters Ten and Eleven study general properties of types of countable sets (the numbers of a second class, according to Cantor). The twelve's chapter is devoted to forming of a scale of alephs, thirteenth chapter studies power of potency. Zhegalkin proved König's theorem for uncountable case, equivalent to axiom of choice in multiplicative form. The monograph ends with a list of paradoxes known by that time. In fact, Zhegalkin made an attempt to build a consistent and complete transfinite number theory. However, he was based on the concept of a finite set without strictly defining it. He also studied numbers that are above class II, it was not in Cantor's works.

## Moscow School of the Theory of Functions and Sets

In 1910, D.F. Egorov launched a workshop in Moscow University in the Theory of Functions. In 1911, the history of Moscow School of Set Theory started with Egorov's theorem of uniform convergence. This School was headed by Egorov and N.N. Luzin. Luzin's research created a new line, descriptive set theory; research of his students developed numerous lines based on the set theory: the theory of measure, set-theoretic topology, functional analysis, the theory of probability, and many other.

## Petersburg–Odessa. 1914. P.S. Yushkevich

Three basic works of Cantor already translated (not paraphrased) were published in 1914. From 1913 to 1915, Vasiliev was publishing series entitled "New Ideas in Mathematics" in Petrograd. To have Cantor's works translated, he engaged a philosopher and translator of philosophic literature, Pavel Solomonovich Yushkevich (1873–1945), father of Adoph Pavlovich Yushkevich. He translated three of the most characteristic works of Cantor which contained the quintessence of his theory: "Grundlagen einer allgemeinen Mannigfaltigkeitslehre", "Ueber die verschiedenen Standpuncte in bezug auf das Actuelle Unendliche", and "Mitteilungen zur Lehre vom Transfiniten"[Cantor, 1914].

We tell almost nothing here about personal contacts of Russian scientists with Cantor. We would only mention that Cantor was elected foreign member of Kharkov Society of Mathematicians.

Cantor's theory as it had originally appeared (a naïve theory of sets) was revised to make basis for new lines in the theory of functions, theory of measure, functional analysis, set-theoretic topology, and many other branches of mathematics. Certain Russian mathematicians addressed directly fundamentals of the theory of sets. Let us mention a Chuvash mathematician, Isaya Maximovich Maximov (1889–1976) among them. He was a post-graduate student of Luzin, dealt in the theory of sets, theory of numbers, theory of functions, and studied the concept of transfinite space created by him in 1930s.

## Moscow-Novosibirsk. 1968. A.I. Fet. Dramatic fate of the first complete translation of Cantor into Russian.

The story I am about to tell was imparted to me in June 2014 by Liudmila Pavlovna Petrova, widow of A.I. Fet, the first translator of all Cantor's works. Now she lives in Novosibirsk.

Abram Ilyich Fet (1924–2007), mathematician, philosopher, opinion journalist, and brilliant translator, was born in Odessa and graduated from Tomsk University. In 1948, he defended his Candidate Thesis in Moscow. His research advisor was L.A. Lusternik. In 1967, he successfully defended his Ph.D. thesis, which contained the currently known result: Fet's theorem about two geodesics. From 1955, he worked in Novosibirsk. That's what Liudmila Pavlovna told me (the fragments of her letter are published with her consent):

"Whereas you deal with Cantor and history at large, it will probably be interesting for you to know about one episode from history of Cantor's heritage in Russia. A.I. Fet translated not only Cantor's biography written by Frenkel, but all his works. The translation he had done was that of the following publication: Georg Cantor, Ernst Zermelo, ed., Gesammelte Abhandlungen mathematischen und philosophischen inhalts, mit erläuternden anmerkungen sowie mit ergänzungen aus dem briefwechsel Cantor-Dedekind, Berlin, Verlag von Julius Springer, 1932.

This publication included almost all works written by Cantor. Furthermore, there were five letters in the appendix from those Cantor and Dedekind wrote to each other, and Cantor's biography written by Frenkel.

He translated those works in 1969–1970 to earn some money, as in autumn 1968 A.I. Fet was sacked after he signed a letter in defence of illegally convicted and remained unemployed till summer of 1972.

The Contract for the translation was with a Moscow publishing house, "Fizmatlit", in the name of A.V. Gladky, as A.I. was debarred from employment. When the translation was already ready and the publishing house started working on the book, the book was rejected by commission of Pontryagin (not the translation, Cantor's book itself!)."

## L.S. Pontryagin

In 1970, L.S. Pontryagin (1908-1988), academician who made a great contribution in topology and variations calculus, headed a group created by him to form part of a section of the editorial-review board at the Academy of Sciences in the USSR, Chief Editorial Board of Physico-Mathematical Literature at NAUKA Publishing House. This is what he himself wrote: "Already before the group was formed, the section resolved to have G. Cantor's collected works translated into the Russian language. When this resolution was put to vote of the section a second time, this issue got to the group. Before the group started considering it, I.R. Shafarevich met me in the canteen and said: "I do not seem to be a member of the section anymore[2] and therefore would like to warn you regarding the collected works of Cantor. Creation of the theory of sets is unduly assigned to Cantor in whole. In fact, quite a large amount of the work was done by Dedekind. This can be seen in letters Cantor and Dedekind exchanged. Therefore, these letters should be enclosed with Cantor's work."

I started thinking over this suggestion of Shafarevich and concluded that Cantor's works should not be published at all, as it is unreasonable to attract attention of young mathematicians to the theory of sets at the moment.

Very popular in Luzin's times, currently the theory of sets has already lost the edge. The group accepted my suggestion, and the book was rejected. The section readily agreed with us regardless the fact that Cantor's works have already been translated! So we had to pay for the translation services." [Pontryagin, 1998, p. 175].

Liudmila Pavlovna Petrova adds:

---

[2] Shafarevich was expelled from the section as a result of a conflict with Pontryagin, the fact whereof was described by Pontryagin.

"Lev Semenovich was mistaken, the translation services have never been paid for.

The typewritten text of the translated works on 536 pages is kept in our home archives. All formulas, insertions, and colour markings for the publishing house were handwritten by A.I. Fet.

When in 1985 F.A. Medvedev and A.P. Yushkevich[3] translated Cantor's works for NAUKA Publishing House they were not aware of the existence of the already completed translation of Fet (or A.V. Gladky)."

E. N. Savenko wrote about Fet's expertise as a translator as follows:

"The scientist was concerned about the issue of translations his whole life. In 1997, speaking at the conference devoted to this issue, Fet noted that in 1960s, "started the epoch of illiterate translations" [Fet, 1997]. He believed that the reasons for that were in the loss of skills of selection of books to be translated and poor competency of translators, that is to say, their inability to understand the essence of the text being translated caused by poor academic training rather than poor knowledge of the language. A.I. Fet himself, an erudite and a person of keen intellect, possessed unique skills necessary to do quality translations: he would promptly perceive all significant ideas and appropriately lay them down." [Savenko, 2011]

L.P. Petrova added: "He told me that, in his opinion, a good translation of a book in math is a translation which would make this book better. A.I. himself looked upon such translation work as a chance to take a good look at the book he was interested in."

Author's remark: I translated the first biography of Cantor written by A. Fraenkel from the German language. However, when I saw the book translated by A.I. Fet, I was carried away by his lucid and vigorous style that made the text a full-blooded and emotional without distorting the original a single iota. Translators would understand me. I believe that Cantor's works translated by Fet should have been published as well, although we already have a very good translation of 1985 at our disposal.

---

[3] A.P. Yushkevich was an editor, but not a translator.

## Moscow–Leningrad. 1985. F.A. Medvedev

In February 1983 Cantor's works were ready for publication, and in 1985 they were published by NAUKA Publishing House [19]. The publication was prepared by Academician Kolmogorov (1903–1987) and a renowned math historian A.P. Yushkevich (1906–1993). The publication included his basic works in the theory of sets, letters Cantor and Dedekind exchanged, and E. Zermelo's notes to the German publication. The underlying source text was the publication of 1932 edited by Zermelo [Cantor, 1932]. Unlike Zermelo's publication which included five letters of those Cantor and Dedekind wrote to each other, the Russian publication of 1985 includes Cantor's works translated by F.A. Medvedev and 49 letters of the above mathematicians regarding the German publication of E. Noether and G. Cavaillès [Briefwechsel Cantor – Dedekind, 1937].

The Russian publication of 1985 [Cantor, 1985] includes the three Cantor's articles as mentioned above translated by P.S. Yushkevich and published in 1914 in the collection of works entitled "New Ideas in Mathematics"; eleven articles translated by Fedor Andreevich Medvedev, including "Principien einer Theorie der Ordnungstypen. Erste Mitteilung" which was not included in the collection of 1932. It was found by A. Grattan-Guinness as a manuscript kept in Mittag-Leffler Institute in Sweden and published by him in 1970 [Grattan-Guinness, 1970]. This article was written by Cantor in 1884 for *Acta Mathematica*, however, it was rejected by Mittag-Leffler as too philosophic.

Fedor Andreevich Medvedev (1923–1994), mathematician and math historian, author of four books and numerous articles in the history of the theory of sets and work of Cantor himself, devoted his whole life to history of mathematics. Not only did he thoroughly translate Cantor's works, letters he exchanged with Dedekind, and Zermelo's comments, he also added his very valuable notes to Cantor's works. Fedor Andreevich was my teacher; it is thanks to him that I started to study the history of Cantor's theory.

The fate of Russian translations of Cantor's works has lived its 20[th] century history together with Russia. People who touched Cantor's heritage were remarkable, and their names have come down in the history of Russian mathematics.

**Reference in Latin**